\documentclass[runningheads]{llncs}
\usepackage[T1]{fontenc}

\usepackage[colorlinks, linkcolor=blue, citecolor=blue]{hyperref}
\usepackage{graphicx}
\usepackage{xcolor}%
\usepackage{upgreek}
\usepackage{listings}%
\begin{document}
\title{Formalization of the Filter Extension Principle (FEP) in Coq}
%
%
\author{Guowei Dou \and
Wensheng Yu}
\authorrunning{G. Dou \and W. Yu}
%
\institute{Beijing Key Laboratory of Space-ground Interconnection and Convergence,
School of Electronic Engineering, Beijing University of Posts and Telecommunications, Beijing, 100876, China
\email{\{dgw,wsyu\}@bupt.edu.cn}}
\maketitle              
\begin{abstract}
The Filter Extension Principle (FEP) asserts that every filter can be extended to an ultrafilter, which plays a crucial role in the quest for non-principal ultrafilters.
Non-principal ultrafilters find widespread applications in logic, set theory, topology, model theory, and especially non-standard extensions of algebraic structures.
Since non-principal ultrafilters are challenging to construct directly, the Filter Extension Principle, stemming from the Axiom of Choice,
holds significant value in obtaining them. This paper presents the formal verification of the Filter Extension Principle,
implemented using the Coq proof assistant and grounded in axiomatic set theory.
It offers formal descriptions for the concepts related to filter base, filter, ultrafilter and more.
All relevant theorems, propositions, and the Filter Extension Principle itself are rigorously and formally verified.
This work sets the stage for the formalization of non-standard analysis and a specific real number theory.

\end{abstract}
\section{Introduction}

Ultrafilter specifically refers to a type of filter with the maximality property, which is a concept derived from topology.
It has widespread applications in logic, set theory, and model theory, and plays a significant role in non-standard extensions of algebraic structures\cite{CN1974}.

Ultrafilters can be divided into principal and non-principal ultrafilters.
The latter can be utilized to construct the hyper-real numbers that serve as the foundation of non-standard analysis\cite{Rob1974}.
Robinson once quoted G\"odel's statements in the preface of his masterpiece {\it Nonstandard Analysis}\cite{Rob1974}:
``There are good reasons to believe that non-standard analysis, in some version or other, will be the analysis of the future.''
The famous Chinese mathematician, Academician Wu Wenjun, once said: ``Non-standard analysis is the true standard analysis.''\cite{Dau1995}

Besides, Wang noticed a special kind of non-principal ultrafilters in \cite{Wan1998}, named ``non-principal arithmetical ultrafilters''\cite{Wan2000},
which can be used to form an useful non-standard model that has simple construction and superior properties,
and this model can be used to construct real numbers\cite{Wan2018}.
This theory fully leverages the properties of non-principal ultrafilters,
using a non-standard extension approach to bypass the usual rational numbers and directly construct real numbers.

However, constructing non-principal ultrafilters presents a slight challenge.
For each set $A$ and its element $a$, a principle ultrafilter corresponding to $a$ can be directly constructed:
$$F_a = \{u : u \subset A \wedge a \in u\}.$$
While there is no direct construction of non-principal ultrafilters found in various literature.
Generally, their existence requires non-constructive proofs relying on the Axiom of Choice or other assumptions\cite{Tho1973,Wan2016}.
For example, Stanis{\l}aw Marcin Ulam once proved the existence of non-principal ultrafilters over $\omega$ by transfinite induction
according to a well-order, that needs Axiom of Chioce to accomplish\cite{Ula1929,Wan2016}, of the power set of $\omega$.

Another easier method to obtain non-principal ultrafilters lies in the Filter Extension Principle (FEP)\cite{Wan2016,Wan2018},
named the Ultrafilter Theorem in \cite{Tho1973}, which is a consequence of the Axiom of Chioce.
FEP is described as follows in \cite{Wan2018}:

\begin{theorem}[Filter Extension Principle]
If the subset family $G$ of $A$ (i.e., $G \subset 2^A$) possesses ``finite intersection property'':
$$\forall a_1, a_2, \cdots, a_n \in G,\ a_1 \cap a_2 \cap \cdots \cap a_n \neq \emptyset,$$
then there exists an ultrafilter $F$ over $A$ satisfying $G \subset F$.
\end{theorem}

Especially, it can be proven that every filter possesses finite intersection property,
FEP thus asserts that every filter can be extended to an ultrafilter.
Using this principle, non-principal ultrafilters can be conveniently obtained by extending a specific filter, known as the Fr\'{e}chet filter.
The objective of this paper is to introduce the formalization of FEP in the proof assistant Coq.

The formal verification of mathematical theorems has made significant progress in recent years with the development of computer science\cite{Cas2021},
especially with the emergence of proof assistants such as Coq\cite{BC2004}, Isabelle/HOL\cite{NPW2002}, Lean\cite{MU2021}, and others.
Through formal methods, many complex theorems, such as the Four-Color Theorem, Odd Order Theorem and Kepler Conjecture, have already been verified by computers\cite{Gon2008,Gon2013,GAAB2013,HABD2015}.
The Lean4 project initiated by Terence Tao to formalize the proof of the Polynomial Freiman-Rusza Conjecture has also succeeded\cite{Tao2023}.
These achievements contribute to the growing influence of machine-assisted verification of mathematical theorems in the academic community.

According to \cite{Wan2018}, FEP can be divided into three lemmas:

\begin{lemma}
If the subset family $G$ of $A$ (i.e., $G \subset 2^A$) possesses ``finite intersection property'':
$$\forall a_1, a_2, \cdots, a_n \in G,\ a_1 \cap a_2 \cap \cdots \cap a_n \neq \emptyset,$$
then there exists a filter base $B$ over $A$ satisfying $G \subset B$.
\end{lemma}

\begin{lemma}
Every filter base over set $A$ can be extended to a filter.
\end{lemma}

\begin{lemma}
Every filter over set $A$ can be extended to an ultrafilter.
\end{lemma}

In this paper, we will present a rigorous implementation of FEP based on these three lemmas,
serving as theoretical guidelines, each corresponding to a subsection in Section 3.
With FEP, the existence of non-principal ultrafilters can also be straightforwardly verified.

Additionally, the study of filters needs to be conducted within the context of set theories.
In recent years, our team has completed a comprehensive formalization of Morse-Kelley (MK) Axiomatic Set Theory\cite{Kel1955,SY2020,YSF2020} in Coq and is constantly improving it.
MK acknowledges ``classes'', which have a broader scope than sets, as fundamental objects.
In MK, every mathematical object (ordered pair, function, integer, etc.) is a class, and only those classes belonging to some other ones are defined as sets\cite{Kel1955}.
The non-set classes are named ``proper classes''. Therefore, MK is a proper extension of ZFC and is more convenient to utilize\cite{SY2020,YSF2020}.
The work in this paper is grounded in MK formal system.

The paper is organized as follows: Section 2 introduces basic concepts needed later such as filters and ultrafilters, along with their formalization in Coq;
Section 3 presents the entire formal implementation of FEP; Section 4 is conclusion.

\section{Filters}

A filter over set $A$ represents a compatible combination of certain properties of elements in $A$;
an ultrafilter actually represents a maximally compatible combination\cite{Wan2016}.

Marx once stated in {\it Theses on Feuerbach}: ``The essence of man is the sum of all social relations.''
A principal ultrafilter over $A$ is the sum of all properties (relations) that a certain element of $A$ possesses.
This approach to understanding filters parallels the philosophical depth with which people contemplate the essence of humanity,
forging a unity between mathematics and sociology\cite{Wan2016}.

The concept of filters was introduced in \cite{Car1937a,Car1937b} by Henri Paul Cartan
in 1937 and subsequently adopted by Bourbaki in their book {\it General Topology}\cite{Bou1995}.

\begin{definition}[Filter]
Assume that $F$ is a family of subsets of $A$ (i.e., $F \subset 2^A$) and satisfies:\\
1) $\emptyset \notin F,\ A \in F,$\\
2) if $a,b \in F$, then $a \cap b \in F$,\\
3) if $a \subset b \subset A$ and $a \in F$, then $b \in F.$\\
$F$ is called a filter over $A$.
\end{definition}

To formalize this definition, two parameters {\tt F} and {\tt A} are required,
and the first condition can be effectively separated into two statements: $\emptyset \notin F$ and $A \in F$.
Besides, the condition that $F$ is a family of subsets of $A$ also cannot be ignored.
Thanks to nice notations in MK formal system\cite{Kel1955,SY2020,YSF2020}, the formalization is extremely close to mathematical language and easy to understand:

\begin{lstlisting}
Definition Filter F A := F $\subset$ pow(A) /\ $\mathrm{\Phi}$ $\notin$ F /\ A $\in$ F
  /\ ($\forall$ a b, a $\in$ F -> b $\in$ F -> (a $\cap$ b) $\in$ F)
  /\ ($\forall$ a b, a $\subset$ b -> b $\subset$ A -> a $\in$ F -> b $\in$ F).
\end{lstlisting}

\noindent where ``{\tt pow(A)}'' represents the power class of $A$.


The notion of ultrafilters, proposed by Frigyes Riesz\cite{Rie1909}, actually predates the
concept of filters and includes one additional condition.

\begin{definition}[Ultrafilter]
A filter $F$ over $A$ is an ultrafilter if it satisfies:
$$\forall a,\ a \subset A\ \Longrightarrow\ a \in F\ \vee\ (A \sim a) \in F,$$
where $A \sim a$ represents $\{u : u \in A \wedge u \notin a\}$.
\end{definition}

Every ultrafilter over possesses maximality, corresponding to a maximally compatible combination of certain properties of elements in $A$.
More specifically, for every ultrafilter $F$, if a filter $G$ contains $F$ (i.e., $F \subset G$) then $G = F$ must hold,
indicating that an ultrafilter cannot be extended into a larger one.
This leads to another equivalent definition for ultrafilters.

\begin{definition}[Maximal Filter]
A filter $F$ over $A$ is a maximal filter if it satisfies:
$$\forall G,\ G\ \mbox{is a filter over}\ A,\ F \subset G\ \Longrightarrow\ G = F.$$
\end{definition}

It is formally verified in Coq that the definition of ultrafilter and maximal filter are totally equivalent:

\begin{lstlisting}
Definition ultraFilter F A := Filter F A
  /\ ($\forall$ a, a $\subset$ A -> a $\in$ F \/ (A $\sim$ a) $\in$ F).
Definition maxFilter F A := Filter F A
  /\ ($\forall$ G, Filter G A -> F $\subset$ G -> G = F).

Corollary ultraFilter_Equ_maxFilter :
  $\forall$ F A, ultraFilter F A <-> maxFilter F A.
\end{lstlisting}

Ultrafilters can be classified into principal ultrafilters and non-principal ones.
The principal ultrafilters are straightforward to construct.

\begin{definition}[Principal Ultrafilter]
For every $a \in A$, the following set
$$\{u : u \subset A \wedge a \in u\},$$
denoted as $F_a$, is an ultrafilter. Each $F_a$, corresponded to the element $a$ of $A$, is called a principal ultrafilter over $A$.
\end{definition}

$F_a$ consists of all subsets of $A$ that includes the element $a$, it can be regarded as the sum of all properties of $a$ in $A$.
Formalizing this definition requires two steps: the first is to construct the set $F_a$, and the second is to verify that $F_a$ is indeed an ultrafilter.

\begin{lstlisting}
Definition F A a := \{ $\uplambda$ u, u $\subset$ A /\ a $\in$ u \}.
\end{lstlisting}

The construction of $F_a$ requires two parameters {\tt A} and {\tt a}.
Mathematically, $a$ should be a member of $A$ ($a \in A$) and $A$ should be a set,
but in the formal definition these information cannot be reflected.
Therefore, when verifying ``{\tt F A a}'' is an ultrafilter, the additional preconditions towards parameters {\tt A} and {\tt a} should be added.

\begin{lstlisting}
Corollary Fn_Corollary2_b : $\forall$ A a, Ensemble A -> a $\in$ A
  -> ultraFilter (F A a) A.
\end{lstlisting}

Non-principal ultrafilters refer to those ultrafilters that are not principal ones,
and they are equivalent to another kind of ultrafilters -- free ultrafilters.

\begin{definition}[Free Ultrafilter]
An ultrafilter $F$ over $A$ is a free ultrafilter if it satisfies:
$$\forall a,\ a\ \mbox{is a finite subset of}\ A\ \Longrightarrow\ a \notin F.$$
\end{definition}

\begin{lstlisting}
Definition free_ultraFilter F A := ultraFilter F A
  /\ ($\forall$ a, a $\subset$ A -> Finite a -> a $\notin$ F).

Corollary free_ultraFilter_Co1 : $\forall$ F A, free_ultraFilter F A
  -> ($\forall$ x, x $\in$ F -> $\sim$ Finite x).
Corollary free_ultraFilter_Co2 : $\forall$ A, Finite A ->
  $\sim$ $\exists$ F, free_ultraFilter F A.
\end{lstlisting}

In formalization code, ``{\tt Finite a}'' comes from MK formal system\cite{Kel1955,SY2020,YSF2020}, indicating that $a$ is a finite set;
``{\tt free\_ultraFilter\_Co1}'' suggests that elements in a free ultrafilter are all infinite sets;
and ``{\tt free\_ultraFilter\_Co2}'' asserts that there exist no free ultrafilters over a finite set.


\begin{definition}[Fr\'{e}chet Filter]
The Fr\'{e}chet filter over set $A$ is denoted as:
$$F_{\sigma} = \{a : a \subset A \wedge A \sim a\ \mbox{is finite} \}.$$
\end{definition}

It can be proven that $F_{\sigma}$ over an infinite set $A$ is just a filter rather than an ultrafilter,
but $F_{\sigma}$ satisfies one of conditions in the definition of free ultrafilters:
$$\forall a,\ a\ \mbox{\it is a finite subset of}\ A\ \Longrightarrow\ a \notin F_{\sigma}.$$

\begin{lstlisting}
Definition F$\sigma$ A := \{ $\uplambda$ a, a $\subset$ A /\ Finite (A $\sim$ a) \}.
Corollary F$\sigma$_is_just_Filter : $\forall$ A, $\sim$ Finite A -> Ensemble A
  -> Filter (F$\sigma$ A) A /\ $\sim$ ultraFilter (F$\sigma$ A) A
    /\ ($\forall$ a, a $\subset$ A -> Finite a -> a $\notin$ (F$\sigma$ A)).
\end{lstlisting}

The formalization of Fr\'{e}chet filter requires one parameter {\tt A},
and similar to the formalization of principal ultrafilters, preconditions towards {\tt A} should be added when describing relevant theorems.

The following proposition reveals the connection between free ultrafilters and Fr\'{e}chet filters.

\begin{proposition}
For an ultrafilter $F$ over infinite set $A$, $F$ is a free ultrafilter if and only if $F_{\sigma}$ is contained in $F$ (i.e., $F_{\sigma} \subset F$).
\end{proposition}

\begin{lstlisting}
Proposition F$\sigma$_and_free_ultrafilter : $\forall$ F A, Ensemble A -> $\sim$ Finite A
  -> ultraFilter F A -> free_ultraFilter F A <-> (F$\sigma$ A) $\subset$ F.
\end{lstlisting}

In the formalization, ``{\tt Ensemble A}'' and ``{\tt $\sim$ Finite A}'' illustrate that $A$ is an infinite set,
they rule the required properties of $A$.

According to this important proposition, utilizing FEP, $F_{\sigma}$ can be extended to an ultrafilter,
which is precisely a free ultrafilter.

\begin{figure}[h]
\centering
\includegraphics[width=0.65\textwidth]{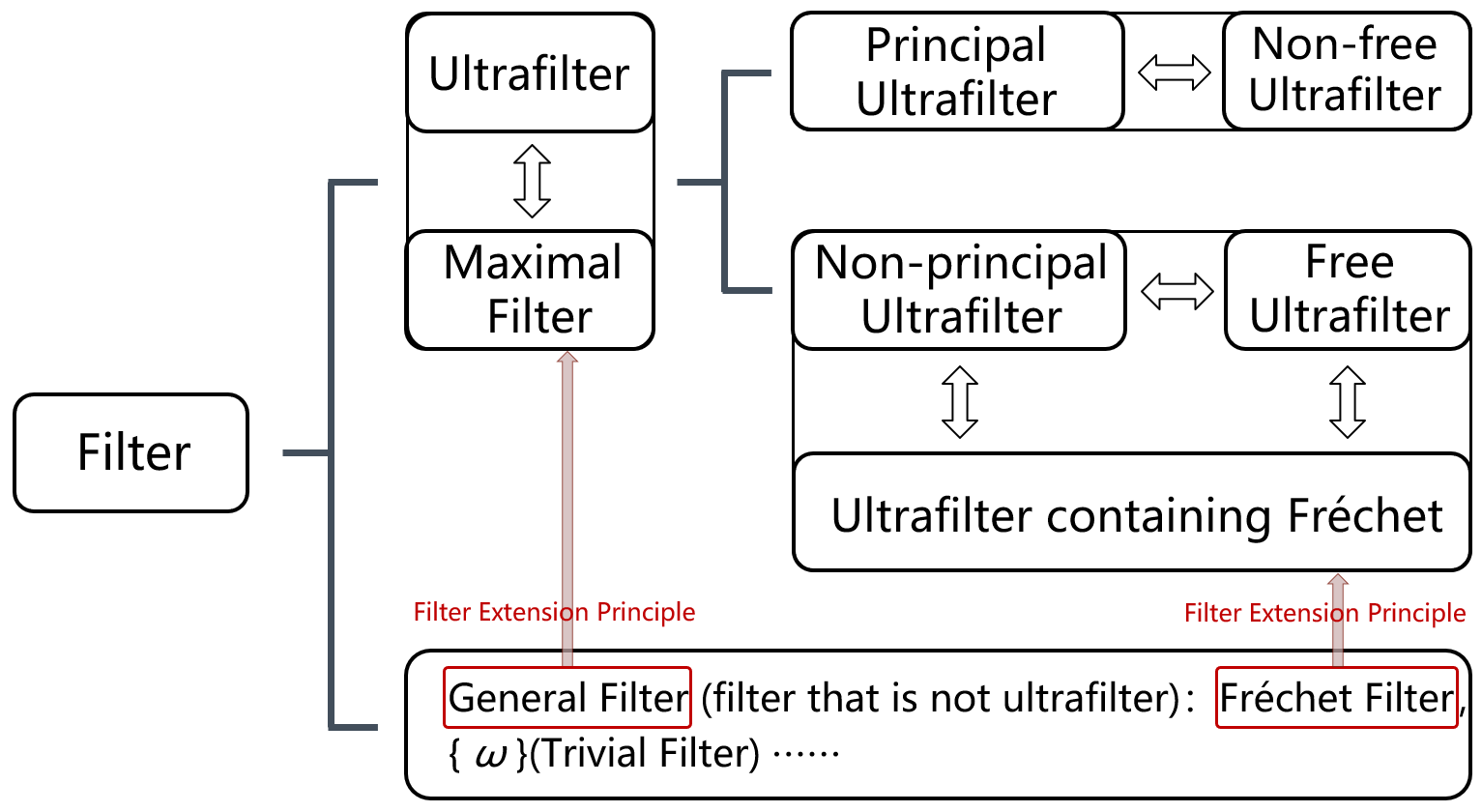}
\caption{relationships among filters}\label{fig1}
\end{figure}

Up to now, we have recognized several concepts related to filters, and Figure.\ref{fig1} presents the relationships among them.
Note that the equivalence between free ultrafilters and non-principal ultrafilters has also been formally verified.
Here we will not go into details; for more information, please refer to the entire code.

Before thoroughly introducing the formalization of FEP,
two more concepts need to be introduced: filter base and finite intersection property.

\begin{definition}[Filter Base]
Assume that $B$ is a non-empty family of subsets of $A$ (i.e., $B \neq \emptyset \wedge B \subset 2^A$) and satisfies:\\
1) $\emptyset \notin B,$\\
2) if $a,b \in B$, then $a \cap b \in B$.\\
$B$ is called a filter base over $A$.
\end{definition}

\begin{lstlisting}
Definition FilterBase B A := B <> $\mathrm{\Phi}$ /\ B $\subset$ pow(A)
  /\ $\mathrm{\Phi}$ $\notin$ B /\ ($\forall$ a b, a $\in$ B -> b $\in$ B -> (a $\cap$ b) $\in$ B).
\end{lstlisting}

\begin{definition}[Finite Intersection Property]
A set $G$ possesses finite intersection property, that is to say, the intersection of a finite number of elements in $G$ is non-empty:
$$\forall a_1, a_2, \cdots, a_n \in G,\ a_1 \cap a_2 \cap \cdots \cap a_n \neq \emptyset.$$
\end{definition}

\begin{lstlisting}
Definition Finite_Intersection G := $\forall$ A, A $\subset$ G -> Finite A -> $\cap$A <> $\mathrm{\Phi}$.
\end{lstlisting}

The formalization for the definition of filter base is straightforward and easily readable;
while describing ``the intersection of a finite number of elements in $G$'' in the definition of finite intersection property requires some skill.
In code, ``{\tt A $\subset$ G}'' and ``{\tt Finite A}'' represent ``a finite number of elements in $G$'';
and ``{\tt $\cap$A}'', whose notation comes from MK formal system\cite{Kel1955,SY2020,YSF2020}, represents the intersection of all these elements,
which is not an empty set (``{\tt $\cap$A <> $\mathrm{\Phi}$}'').

\section{Formalization of the Filter Extension Principle (FEP)}

As mentioned in Introduction, FEP can be divided into three lemmas,
and the idea is to respectively formalize the three lemmas, then verify Theorem 1 (FEP).

\subsection{Formalization of Lemma 1}

The proof of Lemma 1 is a constructive proof. For every subset family $G$ of $A$, if $G$ possesses finite intersection property,
it can be extended to form a filter base in a way that is closed under intersection\cite{Wan2018}:
$$B = \{ a_1 \cap a_2 \cap \cdots \cap a_n : a_1, a_2, \cdots, a_n \in G, n \geq 1 \}.$$
Thus the formalization is divided into two steps: the first is to construct the class $B$, and the second is to prove that $B$ is a filter base containing $G$.

$B$ is formalized with the command ``{\tt\color[RGB]{255,69,0}Definition}'', it requires one parameter ``{\tt G}'' utilized to represent the subset family $G$.

\begin{lstlisting}
Definition FilterBase_from G := \{ $\uplambda$ u, $\exists$ S, S $\subset$ G /\ Finite S /\ u = $\cap$S \}.
\end{lstlisting}


In order to facilitate the use of this formal definition, a notation is introduced:

\begin{lstlisting}
Notation "$\langle$ G $\rangle$$\rightarrow$$^{\mbox{b}}$" := (FilterBase_from G) : filter_scope.
\end{lstlisting}


The next goal is to formally verify that $B$ is indeed a filter base containing $G$.
In formalization, this statement is expressed using the command ``{\tt\color[RGB]{255,69,0}Lemma}'' that needs to be formally proven.

\begin{lstlisting}
Lemma Filter_Extension1 : $\forall$ G A, G <> $\mathrm\Phi$ -> G $\subset$ pow(A)
  -> Finite_Intersection G -> G $\subset$ ($\langle$G$\rangle$$\rightarrow$$^{\mbox{b}}$) /\ FilterBase ($\langle$G$\rangle$$\rightarrow$$^{\mbox{b}}$) A.
\end{lstlisting}

``{\tt G <> $\mathrm\Phi$}'', ``{\tt G $\subset$ pow(A)}'' and ``{\tt Finite\_Intersection G}'' are the preconditions,
which rule that the class $G$ represented by {\tt G} is a subset family of the class $A$ represented by {\tt A}, and $G$ possesses the finite intersection property.
``{\tt G $\subset$ ($\langle$G$\rangle$$\rightarrow$$^{\mbox{b}}$) \verb+/\+ FilterBase ($\langle$G$\rangle$$\rightarrow$$^{\mbox{b}}$) A}''
indicates that $B$, represented by ``{\tt $\langle$G$\rangle$$\rightarrow$$^{\mbox{b}}$}'',
is indeed a filter base over $A$ and is extended from $G$ (i.e., $G \subset B$).

\subsection{Formalization of Lemma 2}

Similar to Lemma 1, the proof of Lemma 2 is also a constructive one.
For every filter base $B$ over set $A$, $B$ can be extended to form a filter in the following way\cite{Wan2018}:
$$F = \{ u : u \subset B\ \wedge\ \exists b, b \in B \wedge b \subset u \}.$$
The formalization is divided into two steps as well.

Firstly, the formal construction of the $F$ requires two parameters ``{\tt B}'' and ``{\tt A}'', respectively utilized to represent the filter base $B$ and the set $A$.

\begin{lstlisting}
Definition Filter_from_FilterBase B A :=
  \{ $\uplambda$ u, u $\subset$ A /\ $\exists$ b, b $\in$ B /\ b $\subset$ u \}.
Notation "$\langle$ B $\mid$ A $\rangle$$^{\mbox{b}}$$\rightarrow$$^{\mbox{f}}$" := (Filter_from_FilterBase B A) : filter_scope.
\end{lstlisting}



Secondly, $F$ is required to be proven a filter over $A$ and is extended from $B$ (i.e., $B \subset F$).

\begin{lstlisting}
Lemma Filter_Extension2 : $\forall$ B A, Ensemble A -> FilterBase B A
  -> B $\subset$ ($\langle$B$\mid$A$\rangle$$^{\mbox{b}}$$\rightarrow$$^{\mbox{f}}$) /\ Filter ($\langle$B$\mid$A$\rangle$$^{\mbox{b}}$$\rightarrow$$^{\mbox{f}}$) A.
\end{lstlisting}

\noindent In the formal description, ``{\tt Ensemble A}'' and ``{\tt FilterBase B A}'' rule that $B$ is a filter base over the set $A$;
``{\tt B $\subset$ ($\langle$B$\mid$A$\rangle$$^{\mbox{b}}$$\rightarrow$$^{\mbox{f}}$) \verb+/\+ Filter ($\langle$B$\mid$A$\rangle$$^{\mbox{b}}$$\rightarrow$$^{\mbox{f}}$) A}'',
as the consequence, states that $F$ (represented by {\tt $\langle$B$\mid$A$\rangle$$^{\mbox{b}}$$\rightarrow$$^{\mbox{f}}$}) is a filter over $A$ and contains the filter base $B$.

According to Lemma 1 and 2, obviously, for every subset family $G$ of $A$, if $G$ possesses finite intersection property,
it can be directly extended to a filter over $A$ in the following way:
$$\langle G \rangle = \{ u : u \subset A\ \wedge\ \exists a_1, a_2, \cdots, a_n \in G, a_1 \cap a_2 \cap \cdots \cap a_n \subset u (n \geq 1) \}.$$

The thought of the formalization for this is inline with that of Lemma 1 and 2:
to formally construct $\langle G \rangle$ firstly then verify it is a filter.

\begin{lstlisting}
Definition Filter_from G A := \{ $\uplambda$ u, u $\subset$ A
  /\ $\exists$ S, S $\subset$ G /\ Finite S /\ $\cap$S $\subset$ u \}.
Notation "$\langle$G$\mid$A$\rangle$$\rightarrow$$^{\mbox{f}}$" := (Filter_from G A) : filter_scope.

Lemma Filter_Extension_1_and_2 : $\forall$ G A, G <> $\mathrm\Phi$ -> G $\subset$ pow(A) -> Ensemble A
  -> Finite_Intersection G -> G $\subset$ ($\langle$G$\mid$A$\rangle$$\rightarrow$$^{\mbox{f}}$) /\ Filter ($\langle$G$\mid$A$\rangle$$\rightarrow$$^{\mbox{f}}$) A.
\end{lstlisting}

\subsection{Formalization of Lemma 3 and FEP}

Different from Lemma 1 and 2, the proof of Lemma 3 requires the Axiom of Choice and is a non-constructive proof.
The formal description for Lemma 3 is straightforward, the only difficulty is how to utilize the Axiom of Choice in the formal proof process.

As is known to us, the Axiom of Choice has many equivalent forms of propositions, and as one of these equivalent statements,
the Zorn's Lemma is utilized to prove Lemma 3 in \cite{Wan2016}.
But in MK formal system, there is no direct description for the Zorn's Lemma.
If complying with the proof line in \cite{Wan2016}, the Zorn's Lemma should be formalized firstly,
and this involves the formalization of some additional definitions and could be a little tedious.

Fortunately, in MK system, another theorem that is equivalent to the Axiom of Choice is included.
It is the Hausdorff Maximal Principle, coded with the name ``{\tt MKT143}'' in MK formal system\cite{Kel1955,SY2020,YSF2020}.
Utilizing ``{\tt MKT143}'', the overall structure of the proof code of Lemma 3 is presented below:

\begin{lstlisting}
Lemma Filter_Extension3 : $\forall$ F A, Filter F A
  -> ($\exists$ F1, F $\subset$ F1 /\ ultraFilter F1 A).
Proof.
  intros. assert (Ensemble F /\ Ensemble A) as []. { $\cdots$ }
  set (M := \{ $\uplambda$ u, F $\subset$ u /\ Filter u A \}).
  assert (Ensemble M). { $\cdots$ }
  $\mbox{\color[RGB]{0,100,0}pose proof}$ H2. apply MKT143 in H3 as [X[[]]].
  assert (F $\in$ X). { $\cdots$ }
  assert (F $\subset$ ($\cup$X)). { $\cdots$ }
  assert (ultraFilter ($\cup$X) A). { $\cdots$ }
  eauto.
Qed.
\end{lstlisting}

\begin{figure}[h]
\centering
\includegraphics[width=0.8\textwidth]{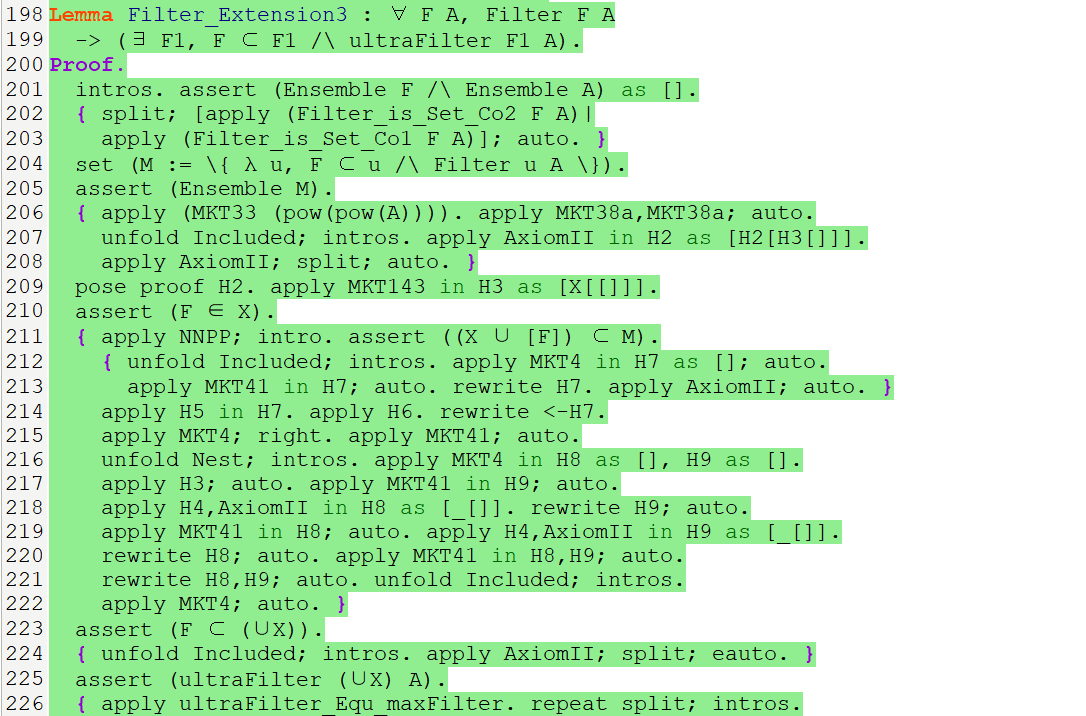}
\caption{partial proof code of Lemma 3}\label{fig2}
\end{figure}

\noindent To highlight the macroscopic process of this proof, the code in ``{\tt \{ $\cdots$ \}}'' is omitted here.
For details, please refer to the complete code. Figure.\ref{fig2} shows the screenshot of partial proof code of Lemma 3.

Now with the use of Lemma 1, 2 and 3, FEP can be formally verified.

\begin{lstlisting}
Theorem Filter_Extension_Principle : $\forall$ G A, G <> $\mathrm\Phi$ -> G $\subset$ pow(A)
  -> Ensemble A -> Finite_Intersection G -> $\exists$ F, G $\subset$ F /\ ultraFilter F A.
Proof.
  intros.
  apply (Filter_Extension_1_and_2 G A) in H2 as []; auto.
  apply Filter_Extension3 in H3 as [F[]]; auto.
  exists F. split; unfold Included; auto.
Qed.
\end{lstlisting}

The formal description of FEP is straightforwardly readable. With the help of the formalization of the previous three lemmas,
the formal proof of FEP is also trivial, it mainly involves the tactic ``{\tt\color[RGB]{0,100,0}apply}'' to call
the formalized lemma ``{\tt Filter\_Extension\_1\_and\_2}'' and ``{\tt Filter\_Extension3}'' that include all the information of the three lemmas.

As an application of FEP, it is also formally verified that there exists a free ultrafilter over each infinite set.
The proof thought has been introduced in Section 2: extending the Fr\'{e}chet filter to an ultrafilter.
Here is the entire formalization of this theorem:

\begin{lstlisting}
Theorem Existence_of_free_ultraFilter : $\forall$ A, Ensemble A
  -> $\sim$ Finite A -> $\exists$ F0, free_ultraFilter F0 A.
Proof.
  intros. $\mbox{\color[RGB]{0,100,0}pose proof}$ H.
  apply F$\sigma$_is_just_Filter in H1 as [H1 _]; auto.
  apply Filter_Extension3 in H1 as [x[]].
  exists x. apply FT2; auto.
Qed.
\end{lstlisting}

In MK system, $\omega$ represents the set of natural numbers\cite{Kel1955,SY2020,YSF2020}, an infinite set,
which satisfies the preconditions in {\tt Existence\_of\_free\_ultraFilter}.
Therefore, there obviously exists a non-principal ultrafilter over natural number set.
Even though this filter is not constructed, it can be instantiated in Coq as:

\begin{lstlisting}
Parameter F0 : Class.
Parameter F0_free_ultraFilter : free_ultraFilter F0 $\upomega$.
\end{lstlisting}

\noindent The first line admits that ``{\tt F0}'' is a constant class in MK system;
and the second line gives ``{\tt F0}" specific mathematical meaning: it is a non-principal ultrafilter over $\omega$.
These two lines of code is actually to admit the existence of a non-principal ultrafilter over $\omega$,
required no formal verification, and the formal theorem {\tt Existence\_of\_free\_ultraFilter} guarantees its consistency.
Then using the constant {\tt F0}, the foundation of non-standard analysis -- hyper-real numbers can be formally constructed\cite{Rob1974}.
The same can be applied to the formalization of non-principal arithmetical ultrafilters mentioned in Introduction,
which serves as the foundation of the formalization of Wang's real number theory\cite{Wan2018}.

\section{Conclusion}

The formalization of the Filter Extension Principle (FEP) comprises four files and involves over two thousand lines of code.
This paper presents the overall architecture of this formalization, as well as the specific formalization process of core definitions and theorems.
Some detailed work is not mentioned in the paper.
For instance, in the proof process of ``$F_{\sigma}$ is a filter'', the proposition that even and odd number sets have the same cardinality is required.
To address this, we have supplemented the formalization related to infinite sets and algebraic operations of natural numbers.
The complete Coq code is available at: {\color{blue}https://github.com/1DGW/Filter-Extension-Principle}

To advance further, the existence of non-principal arithmetical ultrafilters can be formalized based on this work,
which can serve as a preliminary work for the formalization of Wang's real number theory.
Moreover, the formalization of concepts related to filters lays the groundwork for the formalization of non-standard analysis.

\subsubsection{Acknowledgements}

Our work is funded by National Natural Science Foundation (NNSF) of China under Grant 61936008.

%
%
%
%

\end{document}